\documentclass[12pt]{article}

\usepackage{amsmath}
\usepackage{amssymb}
\usepackage{amsfonts}
\usepackage{latexsym}
\usepackage{graphicx}

\newtheorem{proposition}{Proposition}[section]
\newtheorem{theorem}{Theorem}[section]
\newtheorem{corollary}{Corollary}[section]
\newtheorem{remark}{Remark}[section]

\date{ }

\begin{document}

\title{Optimality conditions and Lagrange multipliers for shape and topology optimization problems}

\author{Dan Tiba\thanks{Institute of Mathematics, Romanian Academy, Bucharest, Romania (dan.tiba@imar.ro)}}

\maketitle

\begin{abstract}
  We discuss first order optimality conditions for geometric optimization problems with Neumann boundary conditions and boundary observation.  The methods we develop here are applicable to  large classes of state systems or cost functionals.
	
	Our approach is based on the implicit parametrization theorem and the use of Hamiltonian systems. It establishes equivalence with a constrained optimal control problem and uses Lagrange multipliers under a new simple constraint qualification. In this setting, general functional variations are performed, that combine topological and boundary variations in a natural way.
\end{abstract}



\section{Introduction}
\setcounter{equation}{0}

\label{sec:1}

Let $\mathcal{O}$ be a family of subdomains in $D \subset R^2$, a given bounded domain. We discuss the first order necessary optimality conditions for the following optimal design problem:

\begin{equation}\label{eq:1.1}
  \mathop{\rm Min}\limits_{\Omega \in \mathcal{O}} \mathop{\int}\limits_{\partial\Omega} j(\sigma, y_{\Omega}(\sigma)) d\sigma,
\end{equation} 

\noindent subject to

\begin{equation}\label{eq:1.2}
-\Delta y_{\Omega} + y_{\Omega} = f \quad {\rm in} \; \Omega,
\end{equation}

\begin{equation}\label{eq:1.3}
\frac{\partial y_{\Omega}}{\partial n} = 0 \quad {\rm on} \; \partial\Omega.
\end{equation}

Here, $f \in L^2(D)$ is given and $j(\cdot, \cdot)$ is a Carath\' eodory mapping . More assumptions will be imposed as needed. We work in dimension two since the Poincar\' e-Bendixson theorem \cite{17}, \cite{9} is employed in our periodicity argument for the associated Hamiltonian systems appearing in the description of the unknown geometry. Dimension two is an important case in optimal design. Except this (essential) detail, the proofs are valid in arbitrary dimension.
Moreover, they can be easily extended to other shape optimization problems and our choice for the case (\ref{eq:1.1})-(\ref{eq:1.3}) is motivated by its intrinsic difficulty.

The literature on topology optimization is very rich and we quote just the monographs \cite{1}, \cite{4}, \cite{16} (and their references) devoted to this subject. 
There are many approaches: the topological derivative based on asymptotic analysis of shape functionals in domains with holes, the SIMP method which uses a relaxation procedure, the homogenization approach, level set techniques.
In this paper, we also use level set functions, but no Hamilton-Jacobi equation is needed and no "evolution" of the level sets is taken into account.

Our main tools are general functional variations for the geometry and ordinary differential Hamiltonian systems, via implicit parametrizations of the shapes. In this respect, we refer to \cite{21}, \cite{13}, \cite{15} where such ideas were introduced.
Applications to numerical approaches in optimal design can be found in \cite{22}, \cite{11},  \cite{10}. 

This paper uses as well recent developments from \cite{12}, for instance the differentiability of the period for Hamiltonian systems, with respect to functional variations. However, the key point here is the application of a Lagrange multipliers rule, while in \cite{12} a penalization approximation process plays an essential role. We also quote the paper \cite{2a}, where Lagrange multipliers are obtained under interiority conditions.
Moreover, we are not employing topological asymptotic properties for elliptic equations as in \cite{16}, \cite{8} and we use another topological derivative \cite{15}, \cite{12},  that allows the application of gradient algorithms \cite{10}, \cite{11}.

In this setting, the equivalence of shape optimization problems with certain optimal control problems involving mixed constraints of a special form, is obtained. This may be viewed as a fixed domain approach, but here we have even an equivalence property, not just approximation properties as usual in the literature, \cite{15}, \cite{16}.

Another important point is that our methodology, using functional variations, provides at the computational level simultaneous topological and boundary variations of the unknown geometry.  They are not prescribed, but automatically chosen by the algorithm. As a fixed domain approach it also has clear implementation advantages: the mesh, the mass matrix need not to be updated during the iterations.  Such properties are discussed as well, for instance in \cite{16}, for approximation procedures via different techniques.

It is here that we investigate, for the first time, the impact of general functional variations in the derivation of the necessary optimality conditions in optimal design problems. These necessary conditions have a purely analytic character, that is the geometry is absent from their formulation. We also underline the Frechet differentiability of the operators that we use and the fact that such smooth operations can even generate or close holes, see Rem.\ref{rem:3.2}.

The plan of the paper is as follows. In the next section, we discuss some preliminaries and the optimal control formulation. In Section 3, we obtain abstract optimality conditions via a Lagrange multiplier rule, under a simple constraint qualification condition. The optimality conditions are detailed in the last section. 

Among the advantages of our methodology, we mention its purely analytic character, the generality of the employed variations and the intimate relationship with optimal control theory. We also underline its applicability at the computational level, \cite{12}, \cite{15}. As drawbacks, we indicate the dimension two setting and certain regularity hypotheses, both due to the application of the Poincar\' e-Bendixson theory.

\section{Preliminaries and equivalence}
\label{sec:2}

Let $g \in C(\overline{D})$ be some given function. To it, we associate the open set $\Omega_g \subset D$:

\begin{equation}\label{eq:2.1}
\Omega_g = {\rm int}\{ x \in D; \; g(x) \leq 0 \}
\end{equation}

\noindent that may be not connected and its components may be not simply connected. Notice as well that the set

\begin{equation}\label{eq:2.2}
G = \{ x \in D; \; g(x) = 0 \}
\end{equation}

\noindent may be of positive Lebesgue measure, in this very general framework.

In order to select some connected (not necessarily simply connected) component of $\Omega_g$, we fix a point $\overline{x} \in D$ and define the $ domain$ $\Omega_g$ to be the component that contains $\overline{x}$ in its closure (it may be void). If it is nonvoid, then $g \in C(\overline{D})$ satisfies

\begin{equation}\label{eq:2.3}
g(\overline{x}) \leq 0.
\end{equation}

In fact, in many examples of shape optimization problems, the following supplementary geometric constraint is imposed on the admissible domains $\Omega_g$:

\begin{equation}\nonumber
E \subset \overline{\Omega}_g \Rightarrow g(x) \leq 0 \; {\rm on} \; E,
\end{equation}

\noindent where $E \subset D$ is some given set. Therefore, the condition (\ref{eq:2.3}) and the definition of the domain $\Omega_g$ are natural. 
In the sequel, $\Omega_g$ is always the connected component of the open set defined in (\ref{eq:2.1}), satisfying $\overline{x} \in \overline{\Omega}_g$ and (\ref{eq:2.3}).

It is obvious that more regularity of $\Omega_g$ is necessary in order that (\ref{eq:1.1})-(\ref{eq:1.3}) make sense. We assume  that $g \in C^1(\overline{D})$ and (see (\ref{eq:2.2})):

\begin{equation}\label{eq:2.5}
|\nabla g(x)| > 0, \; \forall x \in G.
\end{equation}

We also impose the condition

\begin{equation}\label{eq:2.6}
g(x) > 0, \; \forall x \in \partial D
\end{equation}

\noindent which ensures that $\partial\Omega_g \cap \partial D = \emptyset$.
We denote by $\mathcal{F} \subset C^1(\overline{D})$ the family of functions satisfying (\ref{eq:2.3}), (\ref{eq:2.5}), (\ref{eq:2.6}). It is a cone in $C^1(\overline{D})$.

By the implicit functions theorem, we get that $\partial\Omega_g$ is of class $C^1$ and

\begin{equation}\label{eq:2.7}
\Omega_g = \{ x \in D, \; g(x) < 0 \},
\end{equation}

\noindent $\partial\Omega_g = G$ and has null Lebesgue measure. Relation (\ref{eq:2.7}) assumes that $g>0$ outside the domain $\Omega_g$, that can be obtained by adding to $g$ the squared distance function to $\Omega_g$ multiplied by a convenient constant. Obviously, an infinity of $g \in \mathcal{F}$ yield the same $\Omega_g$.

If more regularity is imposed on $\mathcal{F}$, then $\partial\Omega_g$ becomes more regular. In general, $\partial D$ is assumed to have the same regularity. The family $\mathcal{O}$ of all admissible domains for the minimization problem (\ref{eq:1.1})-(\ref{eq:1.3}) is generated by the above procedure starting from the level functions $g \in \mathcal{F}$. These domains are not necessarily simply connected and that's why our approach allows topology optimization in combination with shape optimization. We also underline that $\mathcal{O}$ is a rich family of admissible domains, that is the problem (\ref{eq:1.1})-(\ref{eq:1.3}) is meaningful, see \cite{7}, Ch.2.

Denote by $V_{\varepsilon} = \{ x \in D; \; d(x, G) < \varepsilon \} \subset D$ an $\varepsilon$-neighbourhood of $G$; and $G_{\lambda} = \{ x \in D, \; (g + \lambda h)(x) = 0 \}$, where $\varepsilon > 0$, $h \in \mathcal{F}$, $\lambda \in R$.

\begin{proposition}\label{prop:2.1}
Let $\mathcal{F} \subset C(\overline{D})$. There is $\lambda(\varepsilon) > 0$ such that, for $\lambda \in R$, $|\lambda| < \lambda(\varepsilon)$, we have $G_{\lambda} \subset V_{\varepsilon}$.
\end{proposition}

\begin{proposition}\label{prop:2.2}
Let $\mathcal{F} \subset C^2(\overline{D})$. Under conditions (\ref{eq:2.3}), (\ref{eq:2.5}), (\ref{eq:2.6}), $G$ is a finite union of closed curves, without self intersections and disjoint from $\partial D$, globally parametrized by the solution of the Hamiltonian system:

\begin{equation}\label{eq:2.8}
x_1'(t) = -\partial_1 g(x_1(t), x_2(t)), \; t \in I,
\end{equation}

\begin{equation}\label{eq:2.9}
x_2'(t) = \partial_1 g(x_1(t), x_2(t)), \; t \in I,
\end{equation}

\begin{equation}\label{eq:2.10}
x(0) = (x_1(0), x_2(0)) = x_0,
\end{equation}

\noindent where $x_0$ is a point on each component of $G$  and $I$ is an interval around the origin, depending on the respective component. Moreover, the solution of (\ref{eq:2.8})-(\ref{eq:2.10}) is unique.
\end{proposition}

\begin{remark}\label{rem:2.3}
Prop.\ref{prop:2.1} and Prop.\ref{prop:2.2} are very similar to results proved in \cite{22}. A partial (local) extension to arbitrary dimension and the uniqueness result can be found in \cite{21}. In particular, $\mathcal{F}$ is a cone and, for any $g, r \in \mathcal{F}$, we get $g + \lambda r \in \mathcal{F}$ for $\lambda \geq 0$ small, under the conditions of Prop.\ref{prop:2.2} with (\ref{eq:2.3}) replaced by: there is $\overline{x} \in D$ such that

\begin{equation}\label{eq:2.11}
g(\overline{x}) = 0, \; \forall g \in \mathcal{F}.
\end{equation}
\end{remark}

The trajectories of the Hamiltonian system (\ref{eq:2.8})-(\ref{eq:2.10}) are closed, that is periodic and $I$ may be choosen $I = [0, T_g]$, the main period interval, depending on each component of $G$, \cite{22}. Condition (\ref{eq:2.5}) gives the Poincar\' e-Bendixson hypothesis (the absence of equilibrium points on the trajectories) for this system. The Hamiltonian structure has, in fact, the property not to allow the presence of a limit cycle and the solutions are periodic, see \cite{17}, Ch.5,  section 28.

It turns out that the shape optimization problem (\ref{eq:1.1})-(\ref{eq:1.3}) is equivalent with the following constrained optimal control problem, defined in $D$:

\begin{equation}\label{eq:2.12}
\mathop{\rm Min}\limits_{g, u} \mathop{\int}\limits_{0}^{T_g} j(z_g(t), y_g(z_g(t))) |z'_g(t)| dt
\end{equation}

\begin{equation}\label{eq:2.13}
-\Delta y_g + y_g = f + g_{+}^2u \quad {\rm in} \; D,
\end{equation}

\begin{equation}\label{eq:2.14}
y_g = 0 \quad {\rm on} \; \partial D,
\end{equation}

\begin{equation}\label{eq:2.15}
\mathop{\int}\limits_{0}^{T_g} |\nabla y_g(z_g(t)) \cdot \nabla g(z_g(t))|^2 dt = 0.
\end{equation}

In (\ref{eq:2.12})-(\ref{eq:2.15}), $z_g(t) = (z_g^1(t), z_g^2(t))$ denotes the unique solution of the Hamiltonian system (\ref{eq:2.8})-(\ref{eq:2.10}), Constraint (\ref{eq:2.15}) is a simplified equivalent form of (\ref{eq:1.3}):

\begin{equation}\label{eq:2.16}
0 = \mathop{\int}\limits_{\partial\Omega_g} \left|\frac{\partial y_{\Omega_g}}{\partial n}\right|^2 d\sigma = \mathop{\int}\limits_{G} \left|\frac{\partial y_{\Omega_g}}{\partial n}\right|^2 d\sigma = \mathop{\int}\limits_{0}^{T_g} \frac{|\nabla y_g(z_g(t)) \cdot \nabla g(z_g(t))|^2}{|\nabla g(z_g(t))|^2} |z'_g(t)| dt.
\end{equation}

If $G$ has several connected components, then (\ref{eq:2.16}) and (\ref{eq:2.15}) are finite sums (see Prop.\ref{prop:2.2}) of such integrals and an initial condition should be fixed on each component of $G$, for the Hamiltonian system. The same is valid for the cost index (\ref{eq:2.12}) which is an equivalent form of (\ref{eq:1.1}) and may take as well the form of a finite sum if $G = \partial\Omega_g$ has several connected components as it happens in topology optimization.

\begin{theorem}\label{thm:2.4}
For any $g \in \mathcal{F} \subset C^2(\overline{D})$ there is $u_g$ measurable in $D$, not unique, such that $g_+^2 u_g \in L^2(D)$ and the solution of (\ref{eq:2.13}), (\ref{eq:2.14}) coincides in $\Omega_g$ with the solution of (\ref{eq:1.2}), (\ref{eq:1.3}) and satisfies (\ref{eq:2.15}). The associated costs (\ref{eq:1.1}), (\ref{eq:2.12}) coincide as well.

\end{theorem}

\textit{Proof}

Since $\partial\Omega_g$ is in $C^2$ under conditions $\mathcal{F} \subset C^2(\overline{D})$ and (\ref{eq:2.5}), then the unique solution $y$ of (\ref{eq:1.2}), (\ref{eq:1.3}) satisfies $y \in H^2(\Omega_g)$. Let $\widetilde{y} \in H^2(D\setminus \Omega_g)$ be such that $\widetilde{y} = y$, $\frac{\partial \widetilde{y}}{\partial n} = 0$ on $\partial\Omega_g$ and $\widetilde{y} = 0$ on $\partial D$. This is possible due to the trace theorem and $\widetilde{y}$ is not unique.

We denote
$$
u_g = \frac{-\Delta \widetilde{y} - f}{g_+^2} \quad {\rm in} \; D \setminus \overline{\Omega}_g
$$

\noindent and 0 in $\Omega_g$. Here, we may assume $g > 0$ in $D \setminus \overline{\Omega}_g$ without loss of generality.

The concatenation of $y$ and $\widetilde{y}$ satisfies (\ref{eq:2.13}), (\ref{eq:2.14}), (\ref{eq:2.15}) in $D$ and $g_+^2 u_g \in L^2(D)$ by construction.

The equality of the costs (\ref{eq:1.1}), (\ref{eq:2.12}) is obvious. \quad $\Box$

\begin{remark}\label{rem:2.5}
This statement is a variant of the similar results used in other cases in \cite{10}, \cite{11}, \cite{12}. It proves the equivalence of the shape optimization problem (\ref{eq:1.1})-(\ref{eq:1.3}) with the constrained optimal control problem (\ref{eq:2.12})-(\ref{eq:2.15}) defined in $D$. Notice as well that, although in (\ref{eq:1.3}) the Neumann condition is considered, the equivalent formulation in $D$ uses the Dirichlet condition in (\ref{eq:2.14}).
\end{remark}

The existence of some (local) optimal domain $\Omega^{\ast}$ is assumed in this paper, as it is usual in the discussion of optimality conditions, \cite{3}. For general existence results, under minimal geometric regularity assumptions (the segment property), we quote \cite{14}, \cite{19}.

\section{Differentiability and Lagrange multipliers}
\label{sec:3}

In the sequel, we shall use the so called "reduced" problem for the optimal control problem, instead of the formulation given by (\ref{eq:2.12})-(\ref{eq:2.15}) and (\ref{eq:2.8})-(\ref{eq:2.10}). We denote by $A : L^2(D) \to H^2(D) \cap H_0^1(D)$ the isomorphism defined by the Dirichlet problem in $D$ and (\ref{eq:2.13}), (\ref{eq:2.14}) can be written as $y_g = A(f + g_+^2u)$. Then, the constrained optimal control problem can be written as

\begin{equation}\label{eq:3.1}
\mathop{\rm Min}\limits_{g, u}\mathop{\int}\limits_{0}^{T_g} j(z_g(t), A(f + g_+^2u)(z_g(t)))|z'_g(t)| dt,
\end{equation}

\begin{equation}\label{eq:3.2}
\mathop{\int}\limits_{0}^{T_g}\left| \nabla A(f + g_+^2u)(z_g(t)) \cdot \nabla g(z_g(t)) \right|^2 dt = 0.
\end{equation}

The notations $T_g$, $z_g$ are explained in \S 2 and relations (\ref{eq:3.1}), (\ref{eq:3.2}) define the reduced problem.

To fix the ideas and without losing generality, we assume now that $\Omega^{\ast}$ (a local "optimal" domain) is double connected (it has one hole). Then, its boundary $\partial \Omega^{\ast}$ has exactly two components: one "exterior" component and one "interior" component. To further simplify the writing, we assume that the cost functional (\ref{eq:3.1}) is defined just on the "exterior" component. However, the constraint (\ref{eq:3.2}) has to be satisfied on all the components of $\partial \Omega^{\ast}$. We denote by $\overline{z}_g, \widehat{z}_g$ the solutions of the Hamiltonian system (\ref{eq:2.8})-(\ref{eq:2.10}) corresponding to some given initial conditions $\overline{x}, \widehat{x}$, on these two components.

\begin{equation}\label{eq:3.3}
\mathop{\rm Min}\limits_{g, u}\mathop{\int}\limits_{0}^{\widehat{T}_g} j(\widehat{z}_g(t), A(f + g_+^2u)(\widehat{z}_g(t)))|\widehat{z}'_g(t)| dt,
\end{equation}

\begin{equation}\label{eq:3.4}
\mathop{\int}\limits_{0}^{\widehat{T}_g}\left| \nabla A(f + g_+^2u)(\widehat{z}_g(t)) \cdot \nabla g(\widehat{z}_g(t)) \right|^2 dt +\mathop{\int}\limits_{0}^{\overline{T}_g}\left| \nabla A(f + g_+^2u)(\overline{z}_g(t)) \cdot \nabla g(\overline{z}_g(t)) \right|^2 dt = 0.
\end{equation}

In (\ref{eq:3.3}), (\ref{eq:3.4}), the notations $\widehat{T}_g, \overline{T}_g$ are, respectively, the main periods associated to the two Hamiltonian systems used here. All the above considerations from this section are based on Prop.\ref{prop:2.1}, Prop.\ref{prop:2.2},  Thm.\ref{thm:2.4}. In case $\partial \Omega^{\ast}$ has more components, then initial conditions have to be chosen on each of them and the corresponding Hamiltonian systems have to be included in the state system governing the optimal control problem.

In fact, we discuss about the optimal domain $ \Omega^{\ast}$, just for intuition. The problem (\ref{eq:3.3}), (\ref{eq:3.4}) has no reference to the geometry and we denote by $[g^{\ast}, u^{\ast}] \in \mathcal{F} \times L^2(D)$ some (local) optimal pair for (\ref{eq:3.3}) and satisfying (\ref{eq:3.4}). As $\mathcal{F} \subset C^1(\overline{D})$, we get $(g_+^{\ast})^2 \in C(\overline{D})$ and $(g_+^{\ast})^2 u^{\ast}$ is in $L^2(D)$ (the assumption on $u^{\ast}$ is here slightly stronger than in \S 2).

In this analytic setting, we can define functional variations \cite{13}, \cite{15}, around $[g^{\ast}, u^{\ast}]$, $g^{\ast} + \lambda h$, $u^{\ast} + \lambda v$, where $\lambda \in R$, $h \in \mathcal{F}$ and $v \in L^2(D)$. We underline that the condition $u^{\ast} \in L^2(D)$ allows to consider variations $v \in L^2(D)$, independently of $g^{\ast}$, ensuring that the right-hand side in (\ref{eq:2.13}) is in $L^2(D)$.

If $\mathcal{F}$ is given by (\ref{eq:2.5}), (\ref{eq:2.6}) and (\ref{eq:2.11}) (in the form (\ref{eq:3.5}) below), then $g^{\ast} + \lambda h \in \mathcal{F}$ for $|\lambda|$ small, due to the compactness of $\partial D$, $G$ and the Weierstrass theorem. That is, the above variations are "admissible" (the question of the constraint (\ref{eq:3.4}) will be handled via the Lagrange multipliers). We also assume that the above structure (\ref{eq:3.3}), (\ref{eq:3.4}) of the cost and of the constraint (one interior and one exterior component of the boundary) remains the same for the considered perturbed controls $g^{\ast} + \lambda h$, $u^{\ast} + \lambda v$, for $|\lambda|$ small.  This is reasonable due to Prop.\ref{prop:2.1}: the topological structure of $\partial \Omega^{\ast}$ is maintained under the above small perturbations.

We know that $g^{\ast}(\overline{x}) = g^{\ast}(\widehat{x}) = 0$ by the above choice of the initial conditions for (\ref{eq:2.8})-(\ref{eq:2.10}). We take the variations $h \in \mathcal{F}$, around $g^{\ast}$, satisfying as well

\begin{equation}\label{eq:3.5}
h(\overline{x}) = h(\widehat{x}) = 0.
\end{equation}

In fact, (\ref{eq:3.5}) should be understood in the sense that any points $\overline{x}, \widehat{x} \in G^{\ast} = \partial \Omega^{\ast}$ (one on each of the two components of $G^{\ast}$) such that (\ref{eq:3.5}) is fulfilled, may be used. That is the class of variations $h \in \mathcal{F}$ taken into account, is rich.

\begin{remark}\label{rem:3.2}
Such variations are also used in the numerical experiments from \cite{10}, \cite{11}, \cite{12} and the condition (\ref{eq:3.5}) is  easy to  handle at the numerical level, to obtain descent directions. The number of holes can change during the iterations, both by adding new holes or by closing existing holes. Although functional variations are used both here (for optimality conditions) and in computational examples, one should distinguish between the two situations. Here, $\lambda \rightarrow 0$ and we have stability properties for the topology, while in the case of experiments, $\lambda$ is chosen according to the descent property of the cost functional and it is not necessarily "small". Then, topological changes may  appear.
\end{remark}

It is possible, in principle, to remove condition (\ref{eq:3.5}) by defining $\overline{x}_{\lambda} = {\rm proj}_{G_{\lambda}^{\ast}}(\overline{x})$, $\widehat{x}_{\lambda} = {\rm proj}_{G_{\lambda}^{\ast}}(\widehat{x})$ and, automatically, we get $(g + \lambda h)(\overline{x}_{\lambda}) = (g + \lambda h)(\widehat{x}_{\lambda}) = 0$. But, the aim of this section is to prove differentiability properties, for instance of the mapping $z_g(\cdot)$ depending on $g$ (and on the initial conditions $\overline{x}, \widehat{x}, \overline{x}_{\lambda}, \widehat{x}_{\lambda}$). It is clear that differentiability properties of the ${\rm proj}_{G_{\lambda}^{\ast}}(\cdot)$ or of the corresponding ${\rm dist}(\cdot)$ function have to be valid in this setting. However, as shown, for instance, in Delfour and Zolesio [\cite{7}, p. 169] this may not be true. In order to avoid such technicalities, we have imposed (\ref{eq:3.5}), which still ensures a large class of variations.

We denote by $x_{\lambda} = (x_{1\lambda}, x_{2\lambda})$ the solution of the perturbed Hamiltonian system

\begin{equation}\label{eq:3.6}
x'_{1\lambda} = -\partial_2 g(x_{1\lambda}, x_{2\lambda}) - \lambda \partial_2 h(x_{1\lambda}, x_{2\lambda}),
\end{equation}

\begin{equation}\label{eq:3.7}
x'_{2\lambda} = \partial_1 g(x_{1\lambda}, x_{2\lambda}) + \lambda \partial_1 h(x_{1\lambda}, x_{2\lambda}),
\end{equation}

\begin{equation}\label{eq:3.8}
(x_{1\lambda}(0), x_{2\lambda}(0)) = x_0.
\end{equation}

In (\ref{eq:3.6})-(\ref{eq:3.8}) and in (\ref{eq:2.8})-(\ref{eq:2.10}) we assume  $x_0$ to be either $\overline{x}$ or $\widehat{x}$ and (\ref{eq:3.5}) is fulfilled. Let $w_{\lambda} = \frac{x_{\lambda} - x}{\lambda}$.

\begin{proposition}\label{prop:3.3}
Let $g, h \in C^2(\overline{D})$. Then, $x_{\lambda} \to x$ in $C^1(0,T)^2$ and $w_{\lambda} \to w = (w_1, w_2)$ in $C^1(0,T)^2$. The limit $w$ satisfies the system in variations:

\begin{equation}\label{eq:3.9}
w'_1 = -\nabla \partial_2 g(x)\cdot w - \partial_2 h(x) \quad in \; [0, T],
\end{equation}

\begin{equation}\label{eq:3.10}
w'_2 = \nabla \partial_1 g(x)\cdot w + \partial_1 h(x) \quad in \; [0, T],
\end{equation}

\begin{equation}\label{eq:3.11}
w(0) = (w_1(0), w_2(0)) = (0,0).
\end{equation}
\end{proposition}

\textit{Proof} 

Both (\ref{eq:3.6})-(\ref{eq:3.8}) and  (\ref{eq:2.8})-(\ref{eq:2.10}) have periodic solutions, by the above arguments.
Due to this, the stability property $x_{\lambda} \to x$ in $C^1(0,T)^2$, for some arbitrary $T > 0$ is standard, Barbu \cite{2}, Pontryagin \cite{17}. Subtracting (\ref{eq:3.6})-(\ref{eq:3.8}) and (\ref{eq:2.8})-(\ref{eq:2.10}) we obtain

\begin{equation}\label{eq:3.12}
w'_{1\lambda} = -\frac{1}{\lambda}[\partial_2 g(x_{\lambda}) - \partial_2 g(x)] - \partial_2 h(x_{\lambda}) \quad {\rm in} \; [0, T],
\end{equation}

\begin{equation}\label{eq:3.13}
w'_{2\lambda} = \frac{1}{\lambda}[\partial_1 g(x_{\lambda}) - \partial_1 g(x)] + \partial_1 h(x_{\lambda}) \quad {\rm in} \; [0, T],
\end{equation}

\begin{equation}\label{eq:3.14}
w_{\lambda}(0) = (w_{1\lambda}(0), w_{2\lambda}(0)) = (0,0).
\end{equation}

The mean value theorem allows to replace the parentheses in (\ref{eq:3.12}), (\ref{eq:3.13}), respectively by $\nabla \partial_2 g(\theta_{\lambda})\cdot w_{\lambda}$, $\nabla \partial_1 g(\eta_{\lambda})\cdot w_{\lambda}$, where $\theta_{\lambda}(t), \eta_{\lambda}(t)$ are some points on the segment between $x_{\lambda}(t)$ and $x(t)$ and $\theta_{\lambda}(t) \to x(t)$, $\eta_{\lambda}(t) \to x(t)$ as $\lambda \to 0$.

Moreover, $\nabla \partial_2 g(\theta_{\lambda})$, $\nabla \partial_1 g(\eta_{\lambda})$ have respectively the limits $\nabla \partial_1 g(x)$, $\nabla \partial_2 g(x)$ and $\partial_2 h(x_{\lambda}) \to \partial_2 h(x)$, $\partial_1 h(x_{\lambda}) \to \partial_1 h(x)$ under our assumptions, as $\lambda \to 0$.

From (\ref{eq:3.12})-(\ref{eq:3.14}) and their equivalent formulation using $\theta_{\lambda}, \eta_{\lambda}$, we obtain that $\{ w_{\lambda} \}$ is bounded in $C(0,T)^2$ via the Gronwall lemma. It is also bounded in $C^1(0,T)^2$ due to the above convergences. The Arzela-Ascoli theorem gives $w_{\lambda} \to w$ in $C(0,T)^2$, on a subsequence. One can pass to the limit $w_{\lambda} \to w$ in $C^1(0,T)^2$, in fact, and obtain (\ref{eq:3.9})-(\ref{eq:3.11}). Since the solution of the linear system (\ref{eq:3.9})-(\ref{eq:3.11}) is unique, the limit is valid without taking subsequences. \quad $\Box$

\begin{remark}\label{rem:3.4}
This is a modification of Prop. 6 in \cite{20}. If $g, h \in \mathcal{F}$, then $T > 0$ is arbitrary due to Prop.\ref{prop:2.2}. Moreover, the nonlinear operator $z_g : C^2(\overline{D}) \to C^1(0,T)^2$ is G\^{a}teaux differentiable and its G\^{a}teaux differential is linear bounded, due to (\ref{eq:3.9})-(\ref{eq:3.11}), \cite{20}. It is also continuous with respect to $g \in C^2(\overline{D})$, that is the nonlinear operator $z_g(\cdot)$ is Fr\' echet differentiable in the points of $\mathcal{F} \subset C^2(\overline{D})$.
\end{remark}

\begin{proposition}\label{prop:3.5}
We have:
$$
A(f + g_+^2u) : C^2(\overline{D}) \times L^2(D) \to H^2(D) \cap H_0^1(D)
$$
$$
\nabla A(f + g_+^2u) : C^2(\overline{D}) \times L^2(D) \to H^1(D)^2
$$

\noindent are Fr\' echet differentiable.
\end{proposition}

This is obvious due to the properties of the operator $A$ and of the positive part.
\vspace{2mm}

We denote by $\tilde{C}^2(\overline{D}) = \{g \in C^2(\overline{D});\; g(\hat{x}) = g(\overline{x}) = 0 \}$, a Banach space under the same norm as in $C^2(\overline{D})$. It is clear that $\mathcal{F} \subset \tilde{C}^2(\overline{D})$ is an open cone, due to (\ref{eq:2.5}) and (\ref{eq:2.6}). Notice that (\ref{eq:2.5}) is equivalent with $|g(x)| + |\nabla g(x)| > 0$ in $\overline{D}$ and Weierstrass theorem can be applied here too.

By Prop.\ref{prop:3.3} and Prop.\ref{prop:3.5}, it yields that the nonlinear trace operator $A(f + g_+^2u)(z_g(\cdot)) :[\mathcal{F} \subset \tilde{C}^2(\overline{D})] \times L^2(D) \to H^{3/2}(G)$, given by $g, u \rightarrow y_g|_G$, is Fr\' echet differentiable. Similarly $z'_g :\mathcal{F} \subset \tilde{C}^2(\overline{D}) \to C(0,T)^2$, $\nabla g(z_g(\cdot)) :\mathcal{F} \subset \tilde{C}^2(\overline{D}) \to C^1(G)^2$  (trace of $\nabla g$), $\nabla A(f + g_+^2u)(z_g(\cdot)) :[\mathcal{F} \subset \tilde{C}^2(\overline{D})] \times L^2(D) \to H^{1/2}(G)^2$ (trace of $\nabla y_g$) are Fr\' echet differentiable. In case $G^{\ast} = \partial \Omega^{\ast}$ has more connected components, then in (\ref{eq:3.5}) more initial conditions have to be considered and the definition of $\tilde{C}^2(\overline{D})$ has to be adapted accordingly.

\begin{proposition}[\cite{12}]\label{prop:3.7}
Under condition (\ref{eq:2.5}), the functional $T_g : \mathcal{F} \subset \tilde{C}^2(\overline{D}) \to R$ (the period) is Fr\' echet differentiable.
\end{proposition}

In \cite{12}, the G\^{a}teaux differentiability of the period $T_g$ is proved in Prop.4.2. The Fr\' echet differentiability follows by arguments as in Rem.\ref{rem:3.4}.

We denote by $J : \mathcal{F}\times L^2(D) \subset \tilde{C}^2(\overline{D})\times L^2(D) \to R$ and $S : \mathcal{F}\times L^2(D) \subset \tilde{C}^2(\overline{D})\times L^2(D) \to R$, the integral functionals appearing respectively in (\ref{eq:3.3}), (\ref{eq:3.4}) (the functional $S$ is the sum of the two integral functionals appearing in (\ref{eq:3.4})).

\begin{corollary}\label{cor:3.8}
Assume that $j(\cdot, \cdot) \in C^1(R^2 \times R)$. Then $J, S : \mathcal{F} \times L^2(D) \subset \tilde{C}^2(\overline{D})\times L^2(D) \rightarrow  R$ are Fr\' echet differentiable. 
\end{corollary}

We denote by $S',J'$ their Fr\' echet differentials, respectively. More details and arguments, concerning the above properties, will be discussed in the next section.

We consider now the shape optimization problem (\ref{eq:1.1})-(\ref{eq:1.3}) (in its equivalent control formulation (\ref{eq:2.8})-(\ref{eq:2.10}), (\ref{eq:2.12})-(\ref{eq:2.15})), in the following abstract form:

\begin{equation}\label{eq:3.15}
\mathop{\rm Min}\limits_{g \in \mathcal{F}, u \in L^2(D)} \{ J(g, u); \; S(g, u) = 0 \},
\end{equation}

\noindent where $\mathcal{F} \subset \tilde{C}^2(\overline{D})$ satisfies the assumptions  (\ref{eq:2.5}), (\ref{eq:2.6}) and (\ref{eq:3.5}) follows automatically by the definition of $\tilde{C}^2(\overline{D})$. 

We have slightly restricted the set of admissible auxiliary controls $u \in L^2(D)$ with respect to the assumptions of Thm.\ref{thm:2.4}. But we underline that the class of admissible geometries, described by $\mathcal{F} \subset \tilde{C}^2(\overline{D})$, is not modified. However, the constraint (\ref{eq:2.15}) may restrict the set of admissible controls quite severely. We recall that $[g^{\ast}, u^{\ast}]$ is some (local) optimal control pair, assumed to exist  and we denote by $W^{\ast}$ the closed linear one dimensional subspace in $\tilde{C}^2(\overline{D}) \times L^2(D)$ generated by $[g^{\ast}, u^{\ast}]$.

It is known that a one dimensional subspace in a Banach space admits a complementary subspace. This is due to the Hahn-Banach theorem, \cite{23}. In our case $W^{\ast} = \{ \lambda [g^{\ast}, u^{\ast}]; \; \lambda \in R \}$ and we define $\varphi_0 : W^{\ast} \to R$, $\varphi_0(\lambda[g^{\ast}, u^{\ast}]) = \lambda$ a linear continuous functional. It may be extended to a linear continuous functional $\varphi : \tilde{C}^2(\overline{D}) \times L^2(D) \to R$ such that $\varphi([g^{\ast}, u^{\ast}]) = 1$. Then $P : \tilde{C}^2(\overline{D}) \times L^2(D) \to W^{\ast}$, $P([g,u]) = \varphi([g,u])[g^{\ast}, u^{\ast}]$ is a projection operator and the complementary subspace is $V^{\ast} = (I - P)(\tilde{C}^2(\overline{D}) \times L^2(D))$, i.e. $\tilde{C}^2(\overline{D}) \times L^2(D) = V^{\ast} \oplus W^{\ast}$.

We assume that $S : \tilde{C}^2(\overline{D}) \times L^2(D) = V^{\ast} \oplus W^{\ast} \to R$ has its partial derivative with respect to $W^{\ast}$, $\partial_{W^{\ast}} S(g^{\ast}, u^{\ast}) : W^{\ast} \to R$,  an isomorphism (see \cite{6}, Ch.7).

As $W^{\ast}$ has dimension one, this is equivalent to the directional derivative in $[g^{\ast}, u^{\ast}]$ along $\lambda[g^{\ast}, u^{\ast}]$ to satisfy:

\begin{equation}\label{eq:3.16}
S'(g^{\ast}, u^{\ast})(g^{\ast}, u^{\ast}) \neq 0.
\end{equation}

We have proved:
\begin{theorem}\label{thm:3.10}
Under condition (\ref{eq:3.16}), if $[g^{\ast}, u^{\ast}]$ is a (local) minimizer of (\ref{eq:3.15}), there is $k \in R$, $k \neq 0$, such that
\begin{equation}\label{eq:3.17}
J'(g^{\ast}, u^{\ast}) + kS'(g^{\ast}, u^{\ast}) = 0.
\end{equation}
\end{theorem}

This is the classical Lagrange multipliers approach, under the simple constraint qualification (\ref{eq:3.16}), for equality constrained minimization problems, \cite{6}, Ch.7. Important here is that  $\mathcal{F} \times L^2(D)$ is an open subset of $\tilde{C}^2(\overline{D}) \times L^2(D)$. 

\begin{remark}\label{rem:3.11}
 There are other techniques based on regularization or penalization procedures and involving interiority hypotheses \cite{3}, \cite{5}, \cite{2a} in various forms, but not possible to be applied here. Moreover, the one dimensional character of (\ref{eq:3.16}) provides simplicity, as much as possible in this context. Another advantageous point is that, although obtaining (\ref{eq:3.16}) uses essentially the decomposition $V^{\ast} \oplus W^{\ast}$, the result of (\ref{thm:3.10}) is expressed in terms of the original formulation (\ref{eq:3.15}). The fact that the optimal pair is involved in (\ref{eq:3.16}) is standard in the literature, \cite{3}, \cite{6}.

\end{remark}

In the next section, we shall explain (\ref{eq:3.17}) in full details.

\section{Optimality conditions}
\label{sec:4}

In this section we maintain the hypothesis that $\mathcal{F} \subset \tilde{C}^2(\overline{D})$ and  (\ref{eq:2.5}), (\ref{eq:2.6}), (\ref{eq:3.16}) are satisfied. Moreover, to simplify the setting and to fix the ideas, we assume that the boundary of the optimal domain $\partial\Omega^{\ast} = G^{\ast}$ has just two components. The functional variations $g^{\ast} + \lambda h$,\; $h \in \mathcal{F}$ are satisfying (\ref{eq:3.5}), for some points  $\widehat{x}$, $\overline{x}$ on these components. 

\vspace{1mm}
\noindent
In fact, for any $g \in \mathcal{F}$ and $h \in \tilde{C}^2(\overline{D})$, we have $g + \lambda h \in \mathcal{F}$ for $\lambda$ sufficiently small, depending on $g, h$. The variations $u^{\ast} + \lambda v$ are supposed to satisfy just $v \in L^2(D)$. We denote by $(H)$ all these conditions.

It is a standard result (see \cite{10}) that the corresponding equation in variations associated to (\ref{eq:2.13}), (\ref{eq:2.14}) around $[g^{\ast}, u^{\ast}]$ is 

\begin{proposition}\label{prop:4.1}
The limit of $q_{\lambda} = \frac{1}{\lambda}(y_{[g^{\ast},  u^{\ast}] + \lambda[h,v]} - y_{[g^{\ast}, u^{\ast}]})$ exists in $H^2(D)$ and satisfies

\begin{equation}\label{eq:4.1}
-\Delta q + q = (g_+^{\ast})^2v + 2 g_+^{\ast}u^{\ast}h \quad in \; D,
\end{equation}

\begin{equation}\label{eq:4.2}
q = 0 \quad on \; \partial D.
\end{equation}
\end{proposition}

Prop.\ref{prop:4.1} together with Prop.\ref{prop:3.3} give the system in variations of the state system (\ref{eq:2.8})-(\ref{eq:2.10}), (\ref{eq:2.13}), (\ref{eq:2.14}), on each component of $\partial\Omega^{\ast}$. Due to the periodicity of the Hamiltonian system, the equation (\ref{eq:3.9}), (\ref{eq:3.10}) are valid on any interval $[0, T]$. Following \cite{12}, we also have the formulas for the differentials of the main periods $\widehat{T}_{g^{\ast}}$, $\overline{T}_{g^{\ast}}$ (on each component of $G^{\ast}$), with respect to functional variations, completing Prop.\ref{prop:3.7}. Below, we use the simplified notation ${T}_{g^{\ast}}$. 

\begin{proposition}\label{prop:4.2}
We have

\begin{equation}\label{eq:4.3}
\mathop{\lim}\limits_{\lambda \to 0} \frac{T_{g^{\ast} + \lambda h} - T_{g^{\ast}}}{\lambda} = -\frac{w_2(T_{g^{\ast}})}{(z_{g^{\ast}}^2)'(T_{g^{\ast}})},
\end{equation}

\noindent if $(z_{g^{\ast}}^2)'(T_{g^{\ast}}) \neq 0$.
\end{proposition}

Notice that (\ref{eq:2.5}) yields either this condition or $(z_{g^{\ast}}^1)'(T_{g^{\ast}}) \neq 0$ and in the latter case, relation (\ref{eq:4.3}) has to be replaced by

\begin{equation}\label{eq:4.4}
\mathop{\lim}\limits_{\lambda \to 0} \frac{T_{g^{\ast} + \lambda h} - T_{g^{\ast}}}{\lambda} = -\frac{w_1(T_{g^{\ast}})}{(z_{g^{\ast}}^1)'(T_{g^{\ast}})}.
\end{equation}

We denote by $\theta(g^{\ast}, h)$ this limit, in general.

We give now a detailed formulation of Cor.\ref{cor:3.8}.

\begin{proposition}\label{prop:4.3}
Under the assumptions (H), consider the directional derivative of $S : C^2(\overline{D}) \times L^2(D) \to R$ for perturbations $g^{\ast} + \lambda h$, $u^{\ast} + \lambda v$, $\lambda \in R$, $v \in L^2(D)$, $h \in \tilde{C}^2(\overline{D})$ satisfying $h(\overline{x}) = h(\widehat{x}) = 0$ for some points $\overline{x}, \widehat{x} \in G^{\ast}$ (situated on each of the two components of $G^{\ast}$, corresponding to $g^{\ast}$). It is given by:

\vspace{-3mm}
\begin{multline}\label{eq:4.5}
\hat{\theta}(g^{\ast},h) |\nabla A(f + (g_+^{\ast})^2u^{\ast})(\widehat{x})\cdot\nabla g^{\ast}(\widehat{x})|^2 +
2\mathop{\int}\limits_{0}^{T_{g^{\ast}}} \nabla A(f + (g_+^{\ast})^2u^{\ast})(\widehat{z}_{g^{\ast}}(t))\cdot \nabla g^{\ast}(\widehat{z}_{g^{\ast}}(t))\\
\biggl[ \nabla A(f + (g_+^{\ast})^2u^{\ast})(\widehat{z}_{g^{\ast}}(t))\cdot \nabla h(\widehat{z}_{g^{\ast}}(t)) + \nabla q(\widehat{z}_{g^{\ast}}(t))\cdot \nabla g^{\ast}(\widehat{z}_{g^{\ast}}(t)) +\\
+ H\left[ A(f + (g_+^{\ast})^2u^{\ast})(\widehat{z}_{g^{\ast}}(t)) \right]\hat{w}(t)\cdot \nabla g^{\ast}(\widehat{z}_{g^{\ast}}(t)) +\\
+ \nabla A(f + (g_+^{\ast})^2u^{\ast})(\widehat{z}_{g^{\ast}}(t))\cdot H\left[ g^{\ast}(\widehat{z}_{g^{\ast}}(t)) \right]\hat{w}(t) \biggr]dt +
\end{multline}

\vspace{-10mm}
\begin{multline}\nonumber
\overline{\theta}(g^{\ast},h) |\nabla A(f + (g_+^{\ast})^2u^{\ast})(\overline{x})\cdot\nabla g^{\ast}(\overline{x})|^2 +
2\mathop{\int}\limits_{0}^{T_{g^{\ast}}} \nabla A(f + (g_+^{\ast})^2u^{\ast})(\overline{z}_{g^{\ast}}(t))\cdot \nabla g^{\ast}(\overline{z}_{g^{\ast}}(t))\\
\biggl[ \nabla A(f + (g_+^{\ast})^2u^{\ast})(\overline{z}_{g^{\ast}}(t))\cdot \nabla h(\overline{z}_{g^{\ast}}(t)) + \nabla q(\overline{z}_{g^{\ast}}(t))\cdot \nabla g^{\ast}(\overline{z}_{g^{\ast}}(t)) +\\
+ H\left[ A(f + (g_+^{\ast})^2u^{\ast})(\overline{z}_{g^{\ast}}(t)) \right]\overline{w}(t)\cdot \nabla g^{\ast}(\overline{z}_{g^{\ast}}(t)) +\\
+ \nabla A(f + (g_+^{\ast})^2u^{\ast})(\overline{z}_{g^{\ast}}(t))\cdot H\left[ g^{\ast}(\overline{z}_{g^{\ast}}(t)) \right]\overline{w}(t) \biggr]dt. 
\end{multline}
\end{proposition}
\vspace{10mm}
The notation $H[\cdot]$ is the Hessian matrix and $\hat{\theta}, \overline{\theta}$ are associated to (\ref{eq:4.3}) or (\ref{eq:4.4}) in the points $\hat{x}, \overline{x}$ respectively . Similarly, $\hat{w},  \overline{w}$ refer to the vector solutions of (\ref{eq:3.9})-(\ref{eq:3.11}) associated respectively to the two components of $G^{\ast}$.

\textit{Proof}

It is enough to examine the first term in the definition (\ref{eq:3.4}) of $S$:
$$
\begin{array}{c}
\mathop{\lim}\limits_{\lambda \to 0} \frac{1}{\lambda} \biggl[ \mathop{\int}\limits_{0}^{T_{g^{\ast} + \lambda h}} \left| \nabla A(f + (g^{\ast} + \lambda h)_+^2(u^{\ast} + \lambda v))(\widehat{z}_{g^{\ast} + \lambda h}(t)) \cdot \nabla(g^{\ast} + \lambda h)(\widehat{z}_{g^{\ast} + \lambda h}(t))\right|^2 dt -\\[3mm]
- \mathop{\int}\limits_{0}^{T_{g^{\ast}}} \left| \nabla A(f + (g_+^{\ast})^2u^{\ast})(\widehat{z}_{g^{\ast}}(t)) \cdot \nabla g^{\ast}(\widehat{z}_{g^{\ast}}(t)) \right|^2 dt  \biggr].
\end{array}
$$
Here $\widehat{z}_{g^{\ast}}$, $\widehat{z}_{g^{\ast} + \lambda h}$ denote the solution of the Hamiltonian system (\ref{eq:2.8}), (\ref{eq:2.9}), respectively of the corresponding perturbed Hamiltonian system, with initial condition $\widehat{x} \in G^{\ast} \cap G_{\lambda}^{\ast}$.

We notice first that
$$
\begin{array}{c}
\frac{1}{\lambda} \mathop{\int}\limits_{T_{g^{\ast}}}^{T_{g^{\ast}} + \lambda h} \left| \nabla A(f + (g^{\ast} + \lambda h)_+^2(u^{\ast} + \lambda v))(\widehat{z}_{g^{\ast} + \lambda h}(t)) \cdot \nabla(g^{\ast} + \lambda h)(\widehat{z}_{g^{\ast} + \lambda h}(t))\right|^2 dt \to \\[3mm] 
\to \hat{\theta}(g^{\ast}, h)\left| \nabla A(f + (g_+^{\ast})^2u^{\ast})(\widehat{x}) \cdot \nabla g^{\ast}(\widehat{x}) \right|^2,
\end{array}
$$

\noindent due to Prop.\ref{prop:4.2}, the regularity assumptions and the mean value theorem for integrals.

Concerning the integrals over $[0, T_{g^{\ast}}]$, we have (we indicate one case):
$$
\begin{array}{c}
\frac{1}{\lambda}\biggl\{ \mathop{\int}\limits_{0}^{T_{g^{\ast}}} \biggl| \nabla A(f + (g^{\ast} + \lambda h)_+^2(u^{\ast} + \lambda v))(\widehat{z}_{g^{\ast} + \lambda h}(t)) \cdot \nabla(g^{\ast} + \lambda h)(\widehat{z}_{g^{\ast} + \lambda h}(t)) \biggr|^2 dt - \\[3mm]
- \mathop{\int}\limits_{0}^{T_{g^{\ast}}} \biggl| \nabla A(f + (g_+^{\ast})^2u^{\ast})(\widehat{z}_{g^{\ast}}(t)) \cdot \nabla g^{\ast}(\widehat{z}_{g^{\ast}}(t)) \biggr|^2 dt \biggr\} \to \\[3mm]
\to 2\mathop{\int}\limits_{0}^{T_{g^{\ast}}} \nabla A(f + (g_+^{\ast})^2u^{\ast})(\widehat{z}_{g^{\ast}}(t)) \cdot \nabla g^{\ast}(\widehat{z}_{g^{\ast}}(t)) \biggl[ \nabla A(f + (g_+^{\ast})^2u^{\ast})(\widehat{z}_{g^{\ast}}(t)) \cdot \\[3mm]
\cdot \nabla h(\widehat{z}_{g^{\ast}}(t)) + \nabla q(\widehat{z}_{g^{\ast}}(t)) \cdot \nabla g^{\ast}(\widehat{z}_{g^{\ast}}(t)) + \\[3mm]
+ H\biggl[ A(f + (g_+^{\ast})^2u^{\ast})(\widehat{z}_{g^{\ast}}(t))\biggr]\hat{w}(t) \cdot \nabla g^{\ast}(\widehat{z}_{g^{\ast}}(t)) + \\[3mm]
+ \nabla A(f + (g_+^{\ast})^2u^{\ast})(\widehat{z}_{g^{\ast}}(t))  \cdot H g^{\ast}(\widehat{z}_{g^{\ast}}(t))\hat{w}(t) \biggr] dt.
\end{array}
$$

\noindent This ends the proof. \quad $\Box$

In a similar manner, we establish

\begin{proposition}\label{prop:4.5}
Assume that $j \in C^1(R^3)$ and the conditions (H) hold. Then, the directional derivative of $J$ at $[g^{\ast}, u^{\ast}]$ in the direction $[h,v] \in \mathcal{F} \times L^2(D) \subset  \tilde{C}^2(\overline{D})\times L^2(D)$ is given by

\begin{multline}\label{eq:4.6}
\mathop{\int}\limits_{0}^{T_{g^{\ast}}} \nabla_1 j(\widehat{z}_{g^{\ast}}(t), A(f + (g_+^{\ast})^2u^{\ast})(\widehat{z}_{g^{\ast}}(t))) \cdot \hat{w}(t)\left| \widehat{z}'_{g^{\ast}}(t) \right| dt +\\
+ \mathop{\int}\limits_{0}^{T_{g^{\ast}}} \partial_2 j(\widehat{z}_{g^{\ast}}(t), A(f + (g_+^{\ast})^2u^{\ast})(\widehat{z}_{g^{\ast}}(t))) \biggl[ \nabla A(f + (g_+^{\ast})^2u^{\ast})(\widehat{z}_{g^{\ast}}(t)) \cdot\\
\cdot \hat{w}(t) + q(\widehat{z}_{g^{\ast}}(t))\biggr]\left| \widehat{z}'_{g^{\ast}}(t) \right|dt +\\
+ \mathop{\int}\limits_{0}^{T_{g^{\ast}}} j(\widehat{z}_{g^{\ast}}(t), A(f + (g_+^{\ast})^2u^{\ast})(\widehat{z}_{g^{\ast}}(t))) \frac{\widehat{z}'_{g^{\ast}}(t) \cdot \hat{w}'(t)}{\left| \widehat{z}'_{g^{\ast}}(t) \right|} dt +\\
+ \hat{\theta}(g^{\ast}, h) j(\widehat{x}, A(f + (g_+^{\ast})^2u^{\ast})(\widehat{x})).
\end{multline}
\end{proposition}

In (\ref{eq:4.6}), $\nabla_1 j(\cdot, \cdot)$ denotes the gradient of $j$ with respect to the first two arguments and $\partial_2 j(\cdot, \cdot)$ denotes the derivative of $j$ with respect to the last argument. We have also assumed that the cost $J$ is defined just on the "exterior" component of $G^{\ast}$, that contains $\widehat{x}$. The proof follows the same lines as for (\ref{eq:4.5}). It can be found in \cite{12}, in a more general context.

We can now formulate the adjoint system associated to the constrained control problem (\ref{eq:2.8})-(\ref{eq:2.10}), (\ref{eq:2.12})-(\ref{eq:2.15}). We also use Thm.\ref{thm:3.10}, that is the existence of the Lagrange multiplier $k \in R$ such that (\ref{eq:3.17}) is satisfied.
We first define the adjoint system associated to the elliptic equation (\ref{eq:2.13}), (\ref{eq:2.14}) and to the terms containing $q$ from (\ref{eq:4.5}), (\ref{eq:4.6}):

\begin{multline}\label{eq:4.7}
2k \mathop{\int}\limits_{0}^{T_{g^{\ast}}} \nabla A(f + (g_+^{\ast})^2u^{\ast})(\widehat{z}_{g^{\ast}}(t)) \cdot \nabla g^{\ast}(\widehat{z}_{g^{\ast}}(t)) \nabla q(\widehat{z}_{g^{\ast}}(t)) \cdot \nabla g^{\ast}(\widehat{z}_{g^{\ast}}(t)) dt +\\
+ 2k \mathop{\int}\limits_{0}^{T_{g^{\ast}}} \nabla A(f + (g_+^{\ast})^2u^{\ast})(\overline{z}_{g^{\ast}}(t)) \cdot \nabla g^{\ast}(\overline{z}_{g^{\ast}}(t)) \nabla q(\overline{z}_{g^{\ast}}(t)) \cdot \nabla g^{\ast}(\overline{z}_{g^{\ast}}(t)) dt +\\
+ \mathop{\int}\limits_{0}^{T_{g^{\ast}}} \partial_2 j(\widehat{z}_{g^{\ast}}(t), A(f + (g_+^{\ast})^2u^{\ast}))(\widehat{z}_{g^{\ast}}(t)) q(\widehat{z}_{g^{\ast}}(t))\left| \widehat{z}'_{g^{\ast}}(t) \right| dt =\\
= \mathop{\int}\limits_{D} p(x)(-\Delta q(x) + q(x)) dx,
\end{multline}

\noindent for any $q \in H^2(D) \cap H_0^1(D)$.

Notice that, in (\ref{eq:4.7}), $q$ plays the role of an arbitrary test function and the last term in the left-hand side in (\ref{eq:4.7}) can be rewritten as 

\begin{equation}\label{eq:4.8}
\mathop{\int}\limits_{\widehat{G}^{\ast}} \partial_2 j(\sigma, y^{\ast}(\sigma)) q(\sigma) d\sigma,
\end{equation}

\noindent where $\widehat{G}^{\ast}$ is the component of $G^{\ast}$ containing $\widehat{x}$ and $y^{\ast}$ is the solution of (\ref{eq:2.13}), (\ref{eq:2.14}) associated to $g^{\ast}, u^{\ast}$.
This is also valid for the other terms there, by multiplying and dividing by $|\widehat{z}'_{g^{\ast}}(t)| = |\nabla g^{\ast}(\widehat{z}_{g^{\ast}}(t))|$.

Then, (\ref{eq:4.7}) can be rewritten as:

\begin{multline}\label{eq:4.9}
2k \mathop{\int}\limits_{G^{\ast}} \nabla y^{\ast}(\sigma) \cdot \nabla g^{\ast}(\sigma) \frac{\partial q}{\partial n} d\sigma + \mathop{\int}\limits_{\widehat{G}^{\ast}} \partial_2 j(\sigma, y^{\ast}(\sigma))q(\sigma) d\sigma =\\
= \mathop{\int}\limits_{D} p(x)(-\Delta q(x) + q(x)) dx, \quad \forall q \in H^2(D) \cap H_0^1(D).
\end{multline}

In (\ref{eq:4.9}), we also use that $\frac{\nabla g^{\ast}(\cdot)}{|\nabla g^{\ast}(\cdot)|}$ is the unit normal vector to $G^{\ast}$, that is we replace $\frac{\partial q}{\partial n} = \nabla q(\sigma) \cdot \nabla g^{\ast}(\sigma) \frac{1}{|\nabla g^{\ast}(\sigma)|}$.

The first two terms in (\ref{eq:4.7}) are written as one term since $G^{\ast}$ has two components.

If the right-hand side in (\ref{eq:4.1}) is denoted by $\mu \in L^2(D)$, the correspondence $\mu \to q$ defined by (\ref{eq:4.1}), (\ref{eq:4.2}) is an isomorphism between $L^2(D)$ and $H^2(D) \cap H_0^1(D)$. Consequently, the left-hand side in (\ref{eq:4.9}) is a linear continuous functional of $\mu \in L^2(D)$ and there is a unique $p \in L^2(D)$ such that (\ref{eq:4.9}) is satisfied, by the Riesz theorem.

The function $p \in L^2(D)$ is called a very weak solution (solution by transposition) of the elliptic equation

\begin{equation}\label{eq:4.10}
-\Delta p + p = \xi \quad {\rm in}\; D,
\end{equation}

\noindent where $\xi$ is a functional in the dual of $H^2(D) \cap H_0^1(D)$, expressed by the sum of the boundary integrals in (\ref{eq:4.9}). To (\ref{eq:4.10}), the null condition on $\partial D$ is formally added.

We formulate now the adjoint system that takes into account the terms with $w, w'$,  $\hat{\theta}(g^{\ast}, h)$, $\overline{\theta}(g^{\ast}, h)$ from (\ref{eq:4.5}), (\ref{eq:4.6}). Notice that the term $w'$ may be replaced via (\ref{eq:3.9}), (\ref{eq:3.10}) by quantities including just $w$ and $\nabla h$. In fact, we define  the adjoint systems corresponding to (\ref{eq:3.9})-(\ref{eq:3.11}). It is to be underlined that in (\ref{eq:3.9})-(\ref{eq:3.11}), the state $x$ is to be replaced by $\widehat{z}_{g^{\ast}}$, respectively $\overline{z}_{g^{\ast}}$ (corresponding to the respective component of $G^{\ast}$, taken into account).

Moreover, in relations (\ref{eq:4.5}), (\ref{eq:4.6}), we shall not use the terms containing $\nabla h$ (including the ones appearing from the rewriting of $w'$).

Finally, we write the adjoint system corresponding to the component of $G^{\ast}$ containing $\widehat{x}$. The adjoint state here is 
$\widehat{m} = [\widehat{m}_1, \widehat{m}_2]$:

\begin{multline}\label{eq:4.11}
- \widehat{m}'(t) = \widehat{M}^{\ast}(t)\widehat{m}(t) + \nabla_1 j(\widehat{z}_{g^{\ast}}(t), A(f + (g_+^{\ast})^2u)(\widehat{z}_{g^{\ast}}(t)))\left| \widehat{z}'_{g^{\ast}}(t) \right| +\\
+ \partial_2 j(\widehat{z}_{g^{\ast}}(t), A(f + (g_+^{\ast})^2u^{\ast})(\widehat{z}_{g^{\ast}}(t)))\nabla A(f + (g_+^{\ast})^2u^{\ast})(\widehat{z}_{g^{\ast}}(t))\left| \widehat{z}'_{g^{\ast}}(t) \right| +\\
+ j(\widehat{z}_{g^{\ast}}(t), A(f + (g_+^{\ast})^2u^{\ast}))(\widehat{z}_{g^{\ast}}(t))\widehat{M}^{\ast}(t)\frac{\widehat{z}'_{g^{\ast}}(t)}{\left| \widehat{z}'_{g^{\ast}}(t) \right|} +\\
+ 2k\nabla A(f + (g_+^{\ast})^2u^{\ast})(\widehat{z}_{g^{\ast}}(t))\cdot \nabla g^{\ast}(\widehat{z}_{g^{\ast}}(t))\biggl\{ H^{\ast}\left[ A(f + (g_+^{\ast})^2u^{\ast})(\widehat{z}_{g^{\ast}}(t)) \right]\nabla g^{\ast}(\widehat{z}_{g^{\ast}}(t)) +\\
+ H^{\ast} g^{\ast}(\widehat{z}_{g^{\ast}}(t))\nabla A(f + (g_+^{\ast})^2u^{\ast})(\widehat{z}_{g^{\ast}}(t)) \biggr\},
\end{multline}

\begin{multline}\label{eq:4.12}
\widehat{m}_1(T_{g^{\ast}}) = 0, \; \widehat{m}_2(T_{g^{\ast}}) = -\frac{1}{(\widehat{z}_{g^{\ast}}^2)'(T_{g^{\ast}})}\biggl[ j(\widehat{x}, A(f + (g_+^{\ast})^2u^{\ast}))(\widehat{x}) +\\
+ \biggl|\nabla A(f + (g_+^{\ast})^2u^{\ast})(\widehat{x}) \cdot \nabla g^{\ast}(\widehat{x})\biggr|^2 \biggr].
\end{multline}

In (\ref{eq:4.11}), $\widehat{M}$ is the matrix appearing implicitly in (\ref{eq:3.9}), (\ref{eq:3.10}) for $x(t) = \widehat{z}_{g^{\ast}}(t)$ and $\widehat{M}^{\ast}$ is its adjoint (and $H^{\ast}$ is the adjoint as well).

On the component of $G^{\ast}$ containing $\overline{x}$, the adjoint system is simpler since (\ref{eq:4.6}) is not defined on it. The notations use the natural modification of the notations from (\ref{eq:4.11}), (\ref{eq:4.12}):

\begin{multline}\label{eq:4.13}
-\overline{m}'(t) = \overline{M}^{\ast}(t)\overline{m}(t) + 2k \nabla A(f + (g_+^{\ast})^2u^{\ast})(\overline{z}_{g^{\ast}}(t)) \cdot\\
\cdot \nabla g^{\ast}(\overline{z}_{g^{\ast}}(t))\biggl\{ H^{\ast}\biggl[ A(f + (g_+^{\ast})^2u^{\ast})(\overline{z}_{g^{\ast}}(t)) \biggr] \nabla g^{\ast}(\overline{z}_{g^{\ast}}(t)) +\\
+ H^{\ast}g^{\ast}(\overline{z}_{g^{\ast}}(t))\nabla A(f + (g_+^{\ast})^2u^{\ast})(\overline{z}_{g^{\ast}}(t))\biggr\},
\end{multline} 

\begin{multline}\label{eq:4.14}
\overline{m}_1(T_{g^{\ast}}) = 0, \; \overline{m}_2(T_{g^{\ast}}) = -\frac{1}{(\overline{z}_{g^{\ast}}^2)'(T_{g^{\ast}})}\biggl| \nabla A(f + (g_+^{\ast})^2u^{\ast})(\overline{x}) \cdot \nabla g^{\ast}(\overline{x}) \biggr|^2.
\end{multline}

The adjoint system consists of the equations (\ref{eq:4.9})-(\ref{eq:4.14}). In case $G^{\ast}$ has more components, more equations have to be added in (\ref{eq:3.9}), (\ref{eq:3.10}) and in the adjoint system.

\begin{theorem}\label{thm:4.6}
Under the assumptions (H), there is $k \neq 0$ such that

\begin{multline}\label{eq:4.15}
\mathop{\int}\limits_{0}^{T_{g^{\ast}}} j(\widehat{z}_{g^{\ast}}(t), A(f + (g_+^{\ast})^2u^{\ast}))(\widehat{z}_{g^{\ast}}(t))\frac{\widehat{z}'_{g^{\ast}}(t)}{\left| \widehat{z}'_{g^{\ast}}(t)\right|} \cdot (-\partial_2 h(\widehat{z}_{g^{\ast}}(t)), \partial_1 h(\widehat{z}_{g^{\ast}}(t))) +\\
+ 2k \mathop{\int}\limits_{0}^{T_{g^{\ast}}} \nabla A(f + (g_+^{\ast})^2u^{\ast})(\widehat{z}_{g^{\ast}}(t))\cdot \nabla g^{\ast}(\widehat{z}_{g^{\ast}}(t)) \nabla A(f + (g_+^{\ast})^2u^{\ast})(\widehat{z}_{g^{\ast}}(t))\cdot \nabla h(\widehat{z}_{g^{\ast}}(t)) dt +\\
+ 2k \mathop{\int}\limits_{0}^{T_{g^{\ast}}} \nabla A(f + (g_+^{\ast})^2u^{\ast})(\overline{z}_{g^{\ast}}(t)) \cdot \nabla g^{\ast}(\overline{z}_{g^{\ast}}(t))\nabla A(f + (g_+^{\ast})^2u^{\ast})(\overline{z}_{g^{\ast}}(t))\cdot \nabla h(\overline{z}_{g^{\ast}}(t)) dt +\\
+ \mathop{\int}\limits_{D} p(x)\left[ (g_+^{\ast})^2v + 2g_+^{\ast}u^{\ast}h \right] dx + \mathop{\int}\limits_{0}^{T_{g^{\ast}}} \overline{m}(t) \cdot (-\partial_2 h(\overline{z}_{g^{\ast}}(t)), \partial_1 h(\overline{z}_{g^{\ast}}(t))) dt +\\
+ \mathop{\int}\limits_{0}^{T_{g^{\ast}}} \widehat{m}(t) \cdot (-\partial_2 h(\widehat{z}_{g^{\ast}}(t)), \partial_1 h(\widehat{z}_{g^{\ast}}(t))) dt = 0.
\end{multline}

\noindent for any $v \in L^2(D)$, $h \in \mathcal{F}\subset  \tilde{C}^2(\overline{D})$ such that there are $\widehat{x}, \overline{x} \in G^{\ast}$ situated respectively on each component of $G^{\ast}$, with $h(\widehat{x}) = h(\overline{x}) = 0$.
\end{theorem}

This follows by Thm.\ref{thm:3.10}, Prop.\ref{prop:4.3}, Prop.\ref{prop:4.5} and the definitions of the adjoint systems (\ref{eq:4.9})-(\ref{eq:4.14}) and the systems in variations (\ref{eq:3.9})-(\ref{eq:3.11}), respectively (\ref{eq:4.1}), (\ref{eq:4.2}). The relation (\ref{eq:4.15}) expresses the so-called maximum principle for the problem (\ref{eq:3.3}), (\ref{eq:3.4}), which is a slight modification of the original shape optimization problem (\ref{eq:1.1})-(\ref{eq:1.3}), in the sense that in its equivalent form (\ref{eq:2.12})-(\ref{eq:2.15}) we restrict the class of admissible controls to be $u \in L^2(D)$.

Relation (\ref{eq:4.15}) together with the adjoint systems (\ref{eq:4.9})-(\ref{eq:4.14}) and the state system (\ref{eq:2.13}), (\ref{eq:2.14}) give the optimality conditions for this problem.

\begin{remark}\label{rem:4.7}
We use functional variations and optimal control methods and we don't impose explicitly  classical boundary variations or topological variations \cite{18}, \cite{16}. In the statement of Thm. \ref{thm:4.6}, such variations are intimately combined and the formulation has a purely analytic character, i.e. no geometric condition is involved in it. Moreover, in the recent article  \cite{12}, it is shown that the gradient behind our approach can be effectively used in numerical experiments, including both topology and shape optimization. The corresponding algorithm chooses automatically the type of variation (that is not prescribed) in each iteration  and may perform both topological and boundary variations simultaneously.
\end{remark}

\end{document}